# Multiple Solutions to Exponential Diophantine Equations of Ramanujan-Nagell type: $Cz^2 = D + A.B^n$

Philip Gibbs

Abstract: Ramanujan found five solutions to the exponential Diophantine equation $z^2 + 7 = 2^n$ where $z$ and $n$ are positive integers, and posed the problem of determining whether there are any more. Nagell was the first to prove that there were not. It is natural to then ask if there are other similar Diophantine equations with multiple solutions. In particular, equations of the form $Cz^2 = D + A.B^n$ are known as equations of Ramanujan-Nagell type. Three examples are known with six solutions. Summaries of multiple solutions for different cases are presented. In particular the equation $z^2 = 277665.17^6 + 34^n$ is found to have four solutions and is conjectured to be the only such non-trivial equation of Ramanujan-Nagell type with four solutions when $A = C = 1$ and $B$ is not a power of 2.

## Introduction

In 1913, a few weeks before writing his famous letter to G H Hardy, Srinivasa Ramanujan publicised a problem via the Indian Mathematical Society [1]. It was posed as a short open question:

"$2^n - 7$ is a perfect square for the values 3, 4, 5, 7, 15 of $n$. Find other values if there be any."

35 years passed before Trygve Nagell published the first proof that these are indeed the only solutions [2]. He used the methods of Algebraic Number Theory to reduce the problem over the unique factorisation domain $\mathbb{Z}[\sqrt{-7}]$. The Exponential Diophantine Equation $z^2 + 7 = 2^n$ is now known as the Ramanujan-Nagell equation. Some later simpler proofs using similar methods are also notable [3,4,24].

The Ramanujan-Nagell equation has applications, some of which go beyond Diophantine analysis. Solutions of the equation enumerate examples of Mersenne numbers that are triangular. They also provide the only known two examples of numbers (31 and 8191) that are repunits in two different bases (see Goormaghtigh Conjecture). Furthermore, they are used to shortlist the sizes for perfect binary codes with short Hamming distances.

If interest in triangular numbers is expanded to other polygonal numbers, and Mersenne numbers to any power plus or minus some small number, then the general setting is to look for integer solution to equations of the form,

$$Q(z) = B^n$$

Where $Q(z)$ is a quadratic polynomial. This can be reduced to the form,

$$Cz^2 = D + A.B^n$$

Which is a generalisation of the Ramanujan-Nagell Equation for which $(A, B, C, D) = (1, 2, 1, -7)$.

It is possible to limit interest to the case $C = 1$ without restricting the scope of the problem. Any solutions of $(A, B, C, D)$ are also solutions of $(AC, B, 1, DC)$. If C is square-free the correspondence of solutions is one-to-one. However, it is more practical to retain the variable $C$. For example, some number theorists have studied cases with $A = 1$ but $C \geq 1$.

A broader range of Diophantine equations have sometimes been described as Generalised Ramanujan-Nagell equations. These are sometimes of the form $Q(x) = P_1^{n_1} \ldots P_k^{n_l}$ for some set of small primes $\{P_1, \ldots, P_k\}$ [13]. A variation in a different direction is the Lebesgue-Nagell equation $z^2 + D = y^n$ where $D$ is given and all solutions $(x, y, n)$ are sought. This fits into the more general shape $Cz^k = Dx^m + Ay^n$ allowing $x, y$ and $z$ to vary (encompassing FLT, Beal, Catalan and Pillai), which are properly called generalised Lebesgue-Ramanujan-Nagell equations. For the purposes of this paper the equation $Cz^2 = D + A.B^n$ will be referred to as being of Ramanujan-Nagell type.

Our primary goal is to find examples and families of equations where for fixed $A, B, C$ and $D$, there are multiple solutions varying $(z, n)$. The rational variable $y$ is non-negative to avoid trivial double counting, while $n$ is permitted to be any integer, positive, negative or zero. If we wish to avoid rational solutions where $n$ is negative, $D$ can be multiplied by some square power of $B$ to increase $n$ into a positive only range. The number of solutions with such $z$ and $n$ will be denoted as $K(A, B, C, D)$. This is known to be always finite. A secondary goal might be to find solutions where $y$ and $n$ are large for small coefficients.

There is an extensive literature on Exponential Diophantine equations of Ramanujan-Nagell type. A variety of modern methods can be used to find all solutions for some specific $(A, B, C, D)$ or to find families of equations with higher $K(A, B, C, D)$ for restricted cases. For example, $(1, p, 1, D)$ with $p$ prime is commonly studied.

Methods in use are Algebraic Number Theory, $p$-adics, Diophantine Approximation (e.g. Baker's method [14]), elliptic curves and modular methods [22].

To see that there are only a finite number of solutions, observe that solutions in provide integer points to one of the three elliptic curves (or similar with $A$ and $D$ exchanged to cover negative $n$)

$$y^2 = (AC)x^3 + CD, \quad y^2 = (ABC)x^3 + CD, \quad y^2 = (AB^2C)x^3 + CD$$

The theorem of Seagal implies that an elliptic curve can only have a finite number of integer points [11], which in combination implies that $K(A, B, C, D)$ is finite. In the case of the original Ramanujan-Nagell equation, solving these elliptic curves by standard descent methods provides an alternative route to answering Ramanujan's question [12,5].

It is common when studying such problems to impose further conditions on the parameters $(A, B, C, D)$ to restrict the scope of the problem to cases that are more tractable (e.g it might be required that $gcd(B, D) = 1$ or that $z$ and $n$ must be positive.) Here we will study the more general cases with computational and analytic methods allowing us to uncover new equations with high solution counts. It is not our goal to show that we have all solutions in specific cases.

In particular we use a different application of elliptic curves that enable a better understanding of families of equations with multiple solutions via Mordell 2-descent.

# Primitive Equations

**Definition 1:** A Diophantine equation of Ramanujan-Nagell type is of the form

$$Cz^2 = D + A.B^n$$

Where $A, B, C$ and $D$ are given integers, $A$ and $D$ are non-zero, $C > 0$, $|B| \geq 2$, and we seek solutions in $(z, n)$ where $z$ is rational and $n$ is an integer. The equation is referred to by the quadruple $(A, B, C, D)$. The number of solutions is denoted $K(A, B, C, D)$.

As previously noted, the rational and negative nature of the solutions can be transformed to equivalent integer solutions, but it is convenient to express them in this form.

We will explore cases where $D$ has either positive or negative sign. There is no reason to exclude $A$ or $B$ from having a negative sign in the definition, but these cases are of less interest here. (see [6] for some results where $A < 0$)

To avoid trivially infinite families of Ramanujan-Nagell equations, a definition for primitive equations is useful.

> **Definition 2**: A Diophantine equation of Ramanujan-Nagell type $Cz^2 = D + A.B^n$ is said to be primitive if all of the following hold
>
> (i)   $C$ is positive and square-free
> (ii)  There is no common divisor of $A, C$ and $D$.
> (iii) The greatest common divisor of $A$ and $D$ is square free.
> (iv)  $B^2$ does not divide $D$.
> (v)   $B$ does not divide $A$.
> (vi)  If the square free part of $B$ divides $C$ then $B$ does not divide $D$.

Non-primitive equations can be excluded from lists of cases with multiple solutions because they are represented by smaller primitive equations. The definition takes care not to exclude interesting cases as non-primitive, such as $(A, B, C, D) = (1,2,1,17)$ which has five solutions if negative n is permitted. To avoid negative $n$ this would be written in a non-primitive form $(1,2,1,1088)$ which is sometimes excluded in the literature because $\gcd(D, B) \neq 1$

# Counting solutions: $K(A, B, C, D) \geq 1$

How many equations of Ramanujan-Nagell type $Cz^2 = D + A.B^n$ have at least one solution $(z, n), z \geq 0, n \geq 0$? Our goal is to set an up upper bound on this number for fixed positive $A$ and $B$. This is done separately for the two ranges $0 < D \leq D_M$ or $-D_M \leq D < 0$.

To set up a pattern for higher multiplicities the Ramanujan-Nagell type equation can be written

$$Cz^2 = x + a$$

Where $a = A.B^n$. This is trivially solved as $x = Cz^2 - a$ to provide a value for $D = x$ such that the original equation has at least one solution. These solutions can be counted, at least approximately.

> **Counting Conjecture**: Let $A$, $B$ and $C$ be given positive integers, $A \geq 1, B \geq 2, C \geq 1$, then the number of integer values of $D$ such that the Diophantine equation $z^2 = D + A.B^n$ has at least one integer solution $(z, n)$ in the range $0 < D \leq D_M$ is asymptotically bounded by
>
> $$\|\{D: K(A, B, C, D) > 0, D \leq D_M\}\| < \sqrt{\frac{D_M}{C}} \log_B\left(\frac{4D_M}{A}\right)$$
>
> The number of such solutions in the range $-D_M \leq D < 0$ is bounded by
>
> $$\|\{D: K(A, B, D) > 0, -D_M \leq D < 0\}\| < \frac{\sqrt{D_M}}{\sqrt{C} \ln B} (\pi - 2)$$

> **Proof (incomplete)**:
>
> For positive $D$ we have,

$$\|\{D: K(A,B,C,D) > 0, 0 < D \leq D_M\}\| \leq \sum_{0 < D \leq D_M} K(A,B,C,D)$$

$$= \|\{(z,n): z \geq 0, n \geq 0, 0 \leq Cz^2 - A.B^n \leq D_M\}\|$$

$$\approx \sum_{n=0}^{\infty} \left( \sqrt{\frac{D_M + A.B^n}{C}} - \sqrt{\frac{A.B^n}{C}} \right) < \frac{1}{\sqrt{C}} \int_0^{\infty} \left( \sqrt{D_M + A.B^x} - \sqrt{A.B^x} \right) dx$$

$$= \frac{2}{\sqrt{C} \ln B} \left( \sqrt{D_M} \ln \left( \frac{\sqrt{D_M + A} + \sqrt{D_M}}{\sqrt{A}} \right) - \left( \sqrt{D_M + A} - \sqrt{D_M} \right) \right)$$

$$< \sqrt{\frac{D_M}{C}} \log_B \left( \frac{4D_M}{A} \right)$$

For negative $D$,

$$\|\{D: K(A,B,C,D) > 0, -D_M \leq D < 0\}\| \leq \sum_{-D_M \leq D < 0} K(A,B,C,D)$$

$$= \|\{(z,n): z \geq 0, n \geq 0, 0 \leq A.B^n - Cz^2 \leq D_M\}\|$$

$$\approx \sum_{n=\lceil \log_B(D_M/A) \rceil}^{\infty} \frac{\sqrt{A.B^n}}{C} - \frac{\sqrt{A.B^n - D_M}}{C}$$

$$< \frac{1}{\sqrt{C}} \int_{\log_B(D_M/A)}^{\infty} \left( \sqrt{A.B^x} - \sqrt{A.B^x - D_M} \right) dx = \frac{\sqrt{D_M}}{\sqrt{C} \ln B} (\pi - 2)$$

Although these are derived as approximate upper bounds, the errors are small, and the cases with one solution dominate over those with more than one. The bounds can therefore be taken heuristically as asymptotic counts for large $D_M$.

It can be observed that for almost all equations of Ramanujan-Nagell type as $D$ becomes large, there are no solutions, i.e. $K(A,B,C,D) = 0$.

Furthermore, for fixed $A$, $B$ and $C$ most solutions have positive $D$ when $D_M > \frac{A}{4} e^{\pi-2}$

These bounds include both primitive and non-primitive cases. To count the number of solutions to primitive equations for fixed $A$ and $B$, firstly it is required that $B$ does not divide $A$ otherwise there are no primitive cases. Secondly, we would have to include some multiplicative factors to the bounds.

By counting arguments, we can see that for given $A \geq 1, B \geq 2, C \geq 1$, there are an infinite number of both negative and positive $D$ such that the corresponding Ramanujan-Nagell type equation is primitive and $K(A,B,C,D) = 0$, but how easily can we construct examples?

## Double Solutions: $K(A,B,C,D) \geq 2$

Suppose we are given integers $A \geq 1, B \geq 2, C \geq 1$, and we seek $D$ such that $K(A,B,C,D) \geq 2..$ In other words, $Cz^2 = D + A.B^n$ must have at least two integer solutions $(z,n) \in \{(z_1, n_1), (z_2, n_2)\}, n_2 > n_1$.

$$Cz_1^2 = D + A.B^{n_1}$$

$$Cz_2^2 = D + A.B^{n_2}$$

This can be resolved by Subtracting and factorising but alternatively they can be multiplied to follow a pattern to be used for triple and quadruple solutions. Taking the general case

$$y^2 = (x+a)(x+b)$$

Where $y = Cz_1z_2$, $a = A.B^{n_1}$, $b = A.B^{n_2}$, $x = D$. Completing the square transforms this to

$$y^2 = t^2 + \left(\frac{b-a}{2}\right)^2, \quad t = x + \left(\frac{a+b}{2}\right)$$

$$(y-t)(y+t) = \left(\frac{b-a}{2}\right)^2$$

To solve we factorise

$$b - a = G = CEF$$

Here $E$ and $F$ should the same parity. There are always ways to do this for at least some choices of $n_1$ and $n_2$. For example, $n_1 = 2, n_2 = 4$ always makes $G$ a multiple of 4.

Assign

$$y + t = \frac{CE^2}{2}, \quad y - t = \frac{CF^2}{2}$$

$$t = \frac{CE^2}{4} - \frac{CF^2}{4}$$

$$D = x = t + \left(\frac{a+b}{2}\right) = \frac{CE^2}{4} + \frac{CF^2}{4} - \frac{A.B^{n_1} + A.B^{n_2}}{2}$$

By construction, for given $A \geq 1$ and $B \geq 2$, all values of $D$ such that $K(A, B, C, D) \geq 2$ can be found using this formula. Furthermore, they always belong to a family of solutions

$$D = D(B^n) = \frac{E^2}{4}B^{2n} + \frac{F^2}{4} - \frac{A.B^{n_1} + A.B^{n_2}}{2}B^n$$

Where $D(x)$ is a quadratic polynomial with rational coefficients.

# Triple Solutions: $K(A, B, C, D) \geq 3$

It is well known that families of triple solutions to equations of Ramanujan-Nagell type can be constructed as polynomials of exponentials. See the tables in later sections of this paper for examples and references. Here we will use a novel application of Mordell's results to show that these can be better understood using 2-descent in elliptic curves that factorise over $\mathbb{Q}$. The analysis is distinct from previous elliptic curve methods applied to Ramanujan-Nagell type equations. It uses the following lemma.

> **Mordell's Lemma**: Let $y^2 = (x+a)(x+b)(x+c)$ be an elliptic curve that factorises with distinct integer roots $x = -a, -b, -c$. The abelian Mordell-Weil group $\mathbb{E}$ of this curve is finitely generated with a torsion group containing $\mathbb{Z}_2 \times \mathbb{Z}_2$ formed from the three roots and the identity. Generators of the group can be doubled using the tangent method to form elements of a subgroup $2\mathbb{E}$. Using 2-descent methods it can be shown that the rational points for which $x+a, x+b, x+c$ are each simultaneously rational squares correspond to the elements of $2\mathbb{E}$

This lemma can be applied to equations of Ramanujan-Nagell type as follows. Let $Cz^2 = D + A.B^n$ be such an equation denoted by $(A, B, C, D)$. Then $K(A, B, C, D) \geq 3$ with three solutions $(z_i, n_i), i = 1,2,3$ in non-negative integers if and only if the elliptic curve $y^2 = (x + a)(x + b)(x + c)$ has integer points at $x = CD$ in $2\mathbb{E}$ when $a = CA.B^{n_1}, b = CA.B^{n_2}, c = CA.B^{n_3}$.

This result can be used to construct families of Ramanujan-Nagell type equations with three solutions by identifying rational points and applying the tangent method to construct the doubled points in $2\mathbb{E}$.

**Doubling theorem**: Suppose we have an integer $B$, such that $\lfloor B \rfloor \geq 2$ and distinct integers $n_i, i = 1,2,3$. Let $a = B^{n_1}, b = B^{n_2}, c = B^{n_3}$ be three distinct integers, and suppose $(x, y) = (X, Y)$ is a point on the elliptic curve with $Y \neq 0$.

$$y^2 = (x + a)(x + b)(x + c)$$

where $X$ is an integer. (Note that $Y$ is not necessarily an integer, but $Y^2$ is.)

Then compute integer values for $A, C, D$ using the Mordell doubling formula as follows,

The slope of the tangent $\frac{dy}{dx}$ to the curve at $(X, Y)$ is

$$\frac{dy}{dx} = \lambda = \frac{(X + a)(X + b) + (X + c)(X + b) + (X + b)(X + c)}{2Y}$$

And the doubled point by the tangent method is at

$$x = X_2 = \lambda^2 - (a + b + c + 2X)$$

$X_2$ is a rational number. Define $D$ and $A$ as its numerator and denominator

$$X_2 = \frac{D}{A}$$

Finally let $C$ be the square free part of $A.Y^2$. If the choice of X makes $Y^2$ negative then reverse the sign of $A$ and $D$ to make $C$ positive as required by the definition.

The integers $D + A.B^{n_i} = Cz_i^2$ then all have the same square free part $C$ and it follows that $K(A, B, C, D) \geq 3$.

Notes:

- This procedure can only fail to provide a valid equation of Ramanujan-Nagell type with triple solutions if the doubled point takes the value $X_2 = 0$ since them $D$ will be zero which is excluded by the definition.
- If an equation with $C = 1$ is required then use $(CA, B, 1, CD)$
- If an equation with $A = 1$ is required then X must be chosen so that $\lambda^2$ is an integer.
- A sufficient condition for $A = 1$ or $A = 4$ is that the following three quantities are integers $\alpha = \frac{(X+b)(X+c)}{X+a}, \beta = \frac{(X+a)(X+c)}{X+b}, \gamma = \frac{(x+a)(x+b)}{x+c}$
- 
- If an equation with C=1 and A=1 is required then X must be selected so that Y and $\lambda^2$ are both integers.
- If $P = (X, Y)$ is not a torsion point then the points $2mP$ for natural numbers $m$ provide an infinite sequence of distinct primitive Ramanujan-Nagell equations $(A, B, C, D)$ with triple solutions sharing the same $n_i$ the same $B$ and $C$, and the same square-free part of $A$, but distinct $D$.

**Proof:** This follows from Modell's Lemma.

**Corollary:** Give an integer $B \geq 2$ and three integers $n_i$ where

$$0 \leq n_1 < n_2 < n_3 \leq n_1 + n_2$$

And $n_1 + n_2 + n_3$ is even, then there exist an integers $D > 0$ such that

$D + 4 \cdot B^{n_i} = z_i^2$ is square for $i = 1,2,3$. This generates a Ramanujan-Nagell equation $(4, B, 1, D)$ with at least three solutions $K(4, B, 1, D) \geq 3$.

**Proof:** Apply the doubling theorem with $X = 0$, then $Y = B^{(n_1+n_2+n_3)/2}$. This is an integer so we have $C = 1$.

$$4X_2 = D = B^{(n_1+n_2-n_3)/2} + B^{(n_1+n_3-n_2)/2} + B^{(n_2+n_3-n_1)/2} - 2B^{n_1} - 2B^{n_2} - 2B^{n_2}$$

And $A = 4$. In some cases, $D$ will be a multiple of 4 and we can set $A = 1$.

## Case survey

It is useful to divide the scope of the equations into several cases according to ranges of $A, B, C$ and $D$. For each case we can set a lower limit for $K(A, B, C, D)$ and try to find all solutions using a combination of analysis and computational searches. It turns out to be convenient to divide cases according to the values of $B$ separately from the other coefficients as shown in this table

| Surveyed range of $K(A, B, C, D)$ by case | | | | | | |
|---|---|---|---|---|---|---|
| | $B = 2$ | | $B \geq 3$, odd prime | | $B \geq 4$ composite | |
| | $D < 0$ | $D > 0$ | $D < 0$ | $D > 0$ | $D < 0$ | $D > 0$ |
| $A = 1 \text{ or } 4, C = 1$ | $2 - 5$ | $3 - 5$ | $2 - 3$ | 3 | 3 | 4 |
| $A = 1, C > 1$ | $2 - 3$ | 4 | $3 - 4$ | 4 | 4 | 4 |
| $A > 1, C = 1$ | $4 - 6$ | $5 - 6$ | 4 | 5 | $4 - 5$ | 5 |

Below the cases corresponding to each of these cells is expanded

# Cases ($A = 1$ or $4, C = 1, D < 0$)

The equations with $A = 1$ (sometimes 4) and $C = 1$, and especially with $B = p$ a prime, are the most well studied equations of Ramanujan-Nagell type.

### Subcase $B = 2$

In 1960 Apery showed that for $(1,2,1,D)$ where $D < 0$ is odd, there are at most 2 solutions except for the original case from Ramanujan $(1,2,1,-7)$ with its five solutions [7].

In 1981 Beukers showed that the complete set of primitive cases with 2 solutions $K(1,2,1,D) = 2, D < 0$ odd, is given by $D = -23$ or $D = -(2^k - 1), k > 3$ [8]. Allowing solutions where the square is zero enables the additional case $D = -1$. Allowing $D$ to be even does not appear to allow additional double solution cases, but this is unproven. See also Le [15,16]

### Subcase $B$ an odd prime

The equations $(1, p, 1, D)$ for $p$ an odd prime and $D < 0$, have also been investigated with partial results [7,9,10,13]. Note that Beukers and most others stipulate the extra conditions that $p$ does not divide $D$ and that the square is not zero, but we do not. Le allows $A$ to be 4 which adds some new cases including some triple solutions $(p, D) = (3, -11), (5, -19), (7, -2)$ [17,18]. Allowing D to be a multiple of p does not appear to enable more solutions, but this is unresolved.

## Subcase $B$ composite

Where $B$ is composite there are numerous equations such that $K(1, B, D) = 2, D < 0$, but there appears to be only a handful of equations with three solutions and none with more than 3. In all triple solutions $B$ and $D$ have a common factor. No analysis for composite $B$ is available.

The complete answer can be conjectured from computational searches in the table below.

| | | $A = 1$ or $4, C = 1, D < 0$ | | |
|---|---|---|---|---|
| $A$ | $B$ | $D$ | $K(A, B, 1, D)$ | $n$ |
| | | $B = 2, K(A, B, C, D) \geq 2$ | | |
| 1 | 2 | -1 | 2 | 0,1 |
| 1 | 2 | -7 | 5 | 3,4,5,7,15 |
| 1 | 2 | -23 | 2 | 5,11 |
| 1 | 2 | $-(2^k - 1)$ | 2, $k \geq 4$ | $k, 2k - 2$ |
| | | $B = p, K(A, B, C, D) \geq 2, p$ prime | | |
| 1 | 3 | -2 | 2 | 1,3 |
| 1 | $4k^2 + 1$ | -1 | 2 | 0,1 |
| 1 | $4k^2 + 1$ | $-(3k^2 + 1)$ | 2 $k \geq 1$ | 1,3 |
| 4 | 3 | -11 | 3 | 1,2,5 |
| 4 | 5 | -19 | 3 | 1,2,7 |
| 4 | 7 | -3 | 3 | 0,1,3 |
| 4 | $k^2 + k + 1$ | -3 | 2 $k \geq 1 (\neq 2)$ | 0.1 |
| 4 | $p$ | $-(4p^k - 1)$ | 2 $p \geq 3, k \geq 1$ | $k, 2k$ |
| | | $B, K(A, B, C, D) \geq 3, B$ composite | | |
| 1 | 6 | $-23 \cdot 3^2$ | 3 | 3,4,5 |
| 1 | 6 | $-5183 \cdot 3^2$ | 3 | 6,7,10 |
| 1 | 30 | $-119 \cdot 15^2$ | 3 | 3,4,7 |
| 1 | $9 \cdot 2^{2k} - 2$ | $-(2^{2k+2} - 1)(9 \cdot 2^{2k-1} - 1)^{2k+2}$ | 3 $k \geq 1$ | $2k + 2, 2k + 3, 2k + 5$ |
| 1 | $2^{2k} - 2$ | $-(2^{2k+2} - 9)(2^{2k-1} - 1)^{2k+2}$ | 3 $k \geq 2$ | $2k + 2, 2k + 3, 2k + 5$ |

## Case ($A = 1, C = 1, D > 0$)

### Subcase $B = 2$

In the case where $D > 0$ there is more scope for multiple solutions than for $D < 0$ where only the Ramanujan-Nagell equation itself is exceptional. For the case $B = 2$ it is common for mathematicians to impose the condition that $D$ is odd, and more generally $gcd(B, D) = 0$. This makes analysis easier but it does exclude some interesting cases including the equations with the highest multiplicities. For example, the exceptional equation $z^2 = 1088 + 2^n, n \geq 0$ has 5 solutions, but is often left off lists of solutions because 1088 is even. It has however been noted by Ulas [23]. The challenge should be expanded to find all equations with high multiplicity of solutions $K(A, B, 1, D) \geq 3$ without such constraints.

When $D$ is odd there is a family of equations with 4 solutions and several families with 3 solutions. Again, some of these are excluded from earlier lists because $D$ is even.

### Subcase $B$ an odd prime

The case where $D > 0$ and $B = p$ is prime has been studied by Apery, Beukers and Le [7,9,18]. Le has also studied equations with $A = 4$ which adds more interesting families with triple and quadruple solutions. The usual condition that $p$ does not divide $D$ allows most remaining cases to be analysed by the traditional methods with the conclusion that $K(1, p, 1, D) \leq 3$, and $K(4, p, 1, D) \leq 4$ but it excludes the exceptional case of $z^2 = 117 + 4.3^n$ which has four solutions. A complete analysis without this condition would be welcome.

### Subcase $B$ composite

For this case there are many equations with 3 solutions. It is not obvious whether they can be categorised. However, for $K(A, B, C, D) \geq 4$ there are just a handful of examples. In particular the following equation with 4 solutions $n \in \{6,7,8,10\}$ is exceptional and I have not found any prior mention of it.

$$z^2 = 277665.17^6 + 34^n$$

**Conjecture:** This is the only equation of the form $z^2 = D + B^n$ with at least 4 solutions in non-negative integer pairs $(z, n)$, where $B$ is positive and not a power of 2, and $D$ is not a multiple of $B^2$.

Including $A = 4$ does not appear to add any unrelated equations with quadruple solutions when $B$ is composite.

| \multicolumn{5}{c}{$A = 1\ or\ 4, C = 1, D > 0$} | | | | |
|---|---|---|---|---|
| A | B | D | $K(A,B,1,D)$ | $n$ |
| \multicolumn{5}{c}{$B = 2, K(A,B,C,D) \geq 3$} | | | | |
| 1 | 2 | 1088 | 5 | 0,9,11,12,15 |
| 1 | 2 | 1680 | 4 | 0,8,10,12 |
| 1 | 2 | $2^{2k} - 6.2^k + 1$ | $4\ k \geq 3$ | $3, k+2, k+3, 2k+3$ |
| 1 | 2 | $(2^k - 2^m)^2 - 2.(2^k + 2^m) + 1$ | $3\ k \geq m+2 \geq 4$ | $m+2, k+2, k+m+2$ |
| 1 | 2 | $(2^{4k+2} - 34.2^{2k+1} + 1)/9$ | $3\ k \geq 3$ | $5, 2k+3, 4k+7$ |
| 1 | 2 | $(16.2^{4k+2} - 40.2^{2k+1} + 16)/9$ | $3\ k \geq 3$ | $0, 2k+4, 4k+2$ |
| 1 | 2 | $(2^{3k-1} + 2^{2k-1})^2 + 2^{3k} + 2^{2k}$ | $3\ k \geq 3$ | $0, 4k, 8k-4$ |
| 1 | 2 | $(2^{3k-1} - 2^{2k-1})^2 - 2^{3k} + 2^{2k}$ | $3\ k \geq 2$ | $0, 4k, 8k-4$ |
| 1 | 2 | $2^{4k-2} + 2^{2k}$ | $3\ k \geq 4$ | $0, 3k, 6k-6$ |
| \multicolumn{5}{c}{$B = p, K(A,B,C,D) \geq 3, p$ prime} | | | | |
| 4 | 3 | 117 | 4 | 0,3,4,5 |
| 1 | 3 | $(3^{2k+2} - 10.3^{k+1} + 9)/16$ | $3, k \geq 3$ | $0, k+1, 2k$ |
| 1 | 3 | $(3^{2k} - 14.3^k + 1)/16$ | $3, k \geq 3$ | $1, k, 2k+1$ |
| 1 | $p = 4m^2+1$ | $(p^k - 1)^2/4m^2 - p^k$ | $3, k \geq 2$ | $1, k, 2k+1$ |
| 4 | 3 | $(p^k - p^m) - 2(p^k + p^m) + 1$ | $3, k > m \geq 1$ | $m+2, k+2, k+m+2$ |
| \multicolumn{5}{c}{$B, K(A,B,C,D) \geq 4, B$ composite} | | | | |
| 1 | 8 | 1088 | 4 | 0,3,4,5 |
| 1 | 4 | 1680 | 4 | 0,4,5,6 |
| 1 | 34 | $277665.17^6$ | 4 | 6,7,8,10 |

## Case ($A = 1, C > 1$)

The equation $Cz^2 = D + B^n$ has also been well studied, especially for negative $D$ [19,20,10]. To avoid trivial duplications, we demand that $C$ is square free. In general, the number of solutions is low with a few exceptions such as $3z^2 + 5 = 2^n$ which has three solutions, $2z^2 + 14 = 4^n$ and $3z^2 + 13.5^4 = 10^n$ which have 4 solutions.

For $D > 0$ some quadruple solutions are listed in the tables below.

| C | B | D | $K(1,B,C,D)$ | $n$ |
|---|---|---|---|---|
| | | $A=1, C>1, D<0$ | | |
| | | $B=2, K(A,B,C,D) \geq 3$ | | |
| 7 | 2 | -1 | 3 | 0,3,6 |
| 3 | 2 | -5 | 3 | 3,5,9 |
| | | $B$ prime $K(A,B,C,D) \geq 3$ | | |
| 2 | 3 | -1 | 4 | 0,1,2,5 |
| 6 | 7 | -1 | 3 | 0,1,4 |
| | | $B$ composite $K(A,B,C,D) \geq 4$ | | |
| 2 | 4 | -14 | 4 | 2 3 4 8 |
| 3 | 10 | -8125 | 4 | 4 5 6 7 |

| | | A = 1, C > 1, D > 0 | | |
|---|---|---|---|---|
| C | B | D | $K(1,B,C,D)$ | n |
| | | $B = 2, K(A,B,C,D) \geq 4$ | | |
| 3 | 2 | 11 | 4 | 0,4,6,12 |
| | | B prime $K(A,B,C,D) \geq 4$ | | |
| 2 | 7 | 1 | 4 | −2,0,1,2 |
| | | B composite $K(A,B,C,D) \geq 4$ | | |
| 3 | 4 | 11 | 4 | 0 4 6 12 |
| 2 | 55 | 17 | 4 | 0 1 2 4 |
| 14 | 15 | 671 | 4 | 1 2 3 5 |
| 15 | 10 | 8375 | 4 | 3 4 5 9 |
| 3 | 6 | 1971 | 4 | 3 4 5 9 |
| 87 | 30 | 14975 | 4 | 3 4 6 7 |
| 2 | 6 | 11826 | 4 | 4 5 6 10 |
| 3 | 6 | 29187 | 4 | 3 5 6 7 |
| 6 | 10 | 83750 | 4 | 4 5 6 10 |
| 290 | 30 | 4475250 | 4 | 4 5 7 8 |
| 2 | 6 | 175112 | 4 | 4 6 7 8 |

# Case $(A > 1, C = 1)$

Finally, a more general equation is the form $z^2 = D + A.B^n$ for positive or negative $D$. Here we find the only equations known to have 6 solutions:

$$z^2 + 119 = 2^n \text{ [21]} \qquad z^2 = 7.2^{24} + 57.2^n \text{ [23]} \qquad z^2 = 6601.2^2 + 165.2^n \text{ [23]}$$

Further examples with 5 solutions are shown in the tables. This includes some infinite families with 5 solution where $B = 2$ and $D > 0$. [23].

When $B = UV$ is composite there may be finite families of equivalent solutions of the form

$$z^2 = D.U^{m-2k} + (A.V^{2k}).B^n, 0 \leq k \leq m/2$$

These can all be primitive by our definition, but in the table, we will represent them by the example where $k = 0$. The entries in our list turn out to be equivalent to those of Ulas [23].

| A | B | D | K(A,B,1,D) | n |
|---|---|---|---|---|
| $A > 1, C = 1, D < 0$ | | | | |
| $B = 2, K(A,B,C,D) \geq 5$ | | | | |
| 15 | 2 | -119 | 6 | 3 4 5 6 8 15 |
| $2^{6k} + 1$ | 2 | $2^{6k+3} - 1$ | 5  $k \geq 1$ | 3, 2k+2, 6k+2, 6k+4, 18k+6 |
| $(2^{6k-3} + 1)/9$ | 2 | $(2^{6k} - 1)/9$ | 5  $k \geq 2$ | 3, 2k+1, 6k-1, 6k+1, 18k-3 |
| 11 | 2 | -7 | 5 | 0 3 4 8 9 |
| 35 | 2 | -391 | 5 | 4 5 6 11 14 |
| 117 | 2 | -23 | 5 | -2 4 5 8 13 |
| 1045 | 2 | -4879 | 5 | 3 5 6 7 13 |
| 4097 | 2 | -32767 | 5 | 3 6 14 16 42 |
| $B$ prime $K(A,B,C,D) \geq 4$ | | | | |
| 2 | 3 | -2 | 4 | 0 1 2 5 |
| 20 | 3 | -11 | 4 | 0 1 2 3 |
| 26 | 3 | -218 | 4 | 2 3 7 18 |
| 140 | 3 | -299 | 4 | 1 2 3 6 |
| 460 | 3 | -11 | 4 | -2 1 4 7 |
| 760 | 3 | -71 | 4 | -2 1 3 8 |
| 12 | 5 | -11 | 4 | 0 1 2 6 |
| 924 | 5 | -899 | 4 | 0 1 2 5 |
| $B$ composite $K(A,B,C,D) \geq 5$ | | | | |
| 114 | 4 | -455 | 5 | 1 2 5 6 16 |

| $A > 1, C = 1, D > 0$ | | | | |
|---|---|---|---|---|
| A | B | D | $K(A,B,1,$ | n |
| $B = 2, K(A,B,C,D) \geq 5$ | | | | |
| 57 | 2 | 7 | 6 | -24 -10 -8 -4 0 1 |
| 165 | 2 | 6601 | 6 | -2 3 5 6 8 10 |
| $(2^{6k} - 1)/9$ | 2 | $(2^{6k+3} + 1)/9$ | 5  $k \geq 1$ | 0, 2k + 2, 6k + 2, 6k + 4 18k + 6 |
| 5 | 2 | 41 | 5 | -2 3 4 6 9 |
| 65 | 2 | 14 | 5 | -14 -6 -4 -2 1 |
| 39 | 2 | 217 | 5 | 0 3 4 8 11 |
| 105 | 2 | 1 | 5 | -8 -6 -4 3 4 |
| 185 | 2 | 41 | 5 | -6 -4 3 6 14 |
| 195 | 2 | 1921 | 5 | 0 3 4 9 10 |
| 165 | 2 | 721 | 5 | -4 3 4 6 11 |
| 3641 | 2 | 455 | 5 | -42 -16 -14 -6 0 |
| 4097 | 2 | 1022 | 5 | -32 -12 -10 -4 1 |
| 117 | 2 | 9673 | 5 | -2 3 6 7 18 |
| 1305 | 2 | 2329 | 5 | -18 -8 3 6 10 |
| 2465 | 2 | 161 | 5 | -8 -6 -4 3 4 |
| 1881 | 2 | 7153 | 5 | -12 3 4 9 12 |
| 195 | 2 | 24769 | 5 | 0 4 6 7 9 |
| 273 | 2 | 37417 | 5 | -8 3 7 10 12 |
| 221 | 2 | 40145 | 5 | -2 4 6 8 26 |
| 855 | 2 | 106729 | 5 | 0 3 4 11 14 |
| 2145 | 2 | 108601 | 5 | -4 5 7 8 16 |
| 3705 | 2 | 199881 | 5 | -4 3 4 6 7 |
| 3315 | 2 | 961009 | 5 | 0 4 5 7 10 |
| 11505 | 2 | 8449 | 5 | -22 -4 3 6 7 |
| 16385 | 2 | 40094 | 5 | -38 -14 -12 -2 1 |
| 40145 | 2 | 221 | 5 | -26 -8 -6 -4 2 |
| 55385 | 2 | 119 | 5 | -20 -10 -8 -6 1 |
| $B$ prime $K(A,B,C,D) \geq 5$ | | | | |
| 6 | 5 | 19 | 5 | -4 0 1 2 5 |
| 14 | 5 | 11 | 5 | -4 -2 0 1 2 |
| 28 | 3 | 37 | 5 | -4 -2 1 2 6 |
| 70 | 3 | 46 | 5 | -2 1 2 3 6 |
| 130 | 3 | 94 | 5 | -10 -4 1 5 6 |
| 148 | 3 | 517 | 5 | -4 1 2 5 13 |
| 8740 | 3 | 8749 | 5 | -8 -4 1 7 21 |
| $B$ composite $K(A,B,C,D) \geq 5$ | | | | |
| 7 | 4 | 57 | 5 | 0 2 4 5 12 |
| 455 | 4 | 3641 | 5 | 0 3 7 8 21 |
| 62 | 6 | 208449 | 5 | 2 3 5 6 7 |
| 1513 | 6 | 14953473 | 5 | -4 3 5 6 8 |